**PAPER • OPEN ACCESS**

# Control of a non-stationary dynamic system with estimating a strategy of human resources management by the integral indicators method

View the article online for updates and enhancements.

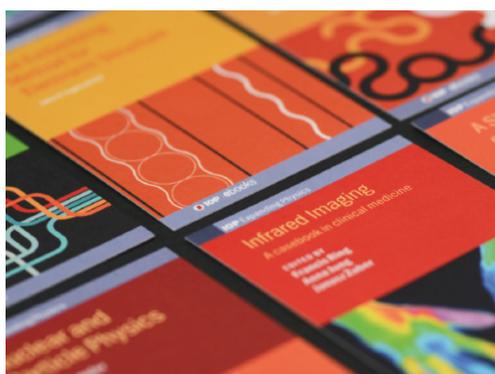

# Control of a non-stationary dynamic system with estimating a strategy of human resources management by the integral indicators method

**S N Masaev[1,4], V V Vingert[1], E V Musiyachenko[1] and Y K Salal[2,3]**

[1] Siberian Federal University, pr. Svobodnyj, 79, Krasnoyarsk, 660041, Russia
[2] South Ural State University, pr. Lenina 76, Chelyabinsk, 454080, Russia
[3] University of Al-Qadisiyah, Diwaniyah, 58002, Iraq

[4] E-mail: faberi@list.ru

**Abstract**. A general idea research is a lack of articles to estimate system indicator of the effectiveness of the strategy of human resources management (HR) at an economic object (enterprise). We are use the method of integral indicators for a comprehensive assessment of the activities of an economic object (enterprise). The economic object is formalized as a non-stationary dynamic system. The system has a dimension of 1.2 million parameters. The parameters of the object under research (business processes) are compared with staff responsibilities. The sanctions mode is set by blocking staff responsibilities in the interval in each time period t. The integral indicator is calculated according to the standard mode (Strategy 1) of the economic object (enterprise) and without blocking the staff responsibilities. Also, the integral indicator is calculated according to the non-standard operating mode of the enterprise (Strategy 2) with the blocking of staff responsibilities. The difference between the integral indicator of Strategy 2 and Strategy 1 is an estimation of the impact of the imposed sanctions. The resource consumption for the restoration of the normal operation of an economic object after the imposed sanctions is give (61.63 million rubles). Example 2 introduces equipment sanctions from America for a new project. An analysis of the indicators of a new project is carried out with the search and use of analog equipment

## 1. Introduction

In control theory, the activity of economic objects, as a system, is researched in various ways: intersectoral balances, vector, parametric and neural network modeling, agent approach, etc. however, the personnel management circuit in such approaches is not identified as a separate object for research. One of the reasons is that the parameters that reveal the behavior of the personnel are not characterized by the normal distribution function, therefore it is difficult to identify and interpret patterns in them. The issue of assessing the state of an economic object as a system is urgent when employees perform their functional duties with various regimes that limit their activity and take into account the influence of environmental parameters.

The management of economic systems involved: V. V. Leontiev and L. V. Kantorovich, A. G. Granberg, A. G. Aganbegyan, V. F. Krotov et al. [1-8]. The research of dynamic system was mainly carried out by foreign authors [9-17]. In 2009, the author proposed integrated indicators [18] for determining the financial crisis of 2008 based on the method of correlation adoptometry. In 2013, integral indicators were used as a separate method, allowing to characterize the activities of the enterprise in various working conditions.

The purpose of this work is: to assess the state of an economic object by an integral indicator as a multidimensional dynamic system, when personnel perform their functional duties in the base mode and sanctions mode with unknown environmental parameters.

To achieve the purpose we perform the tasks:
- Set an economic object (enterprise) as a multidimensional dynamic system (hereinafter the system);
- In the system, assign to the functions to be performed a list of the corresponding functional duties of the personnel from the job descriptions, thereby forming a set of functional responsibilities characterizing the basic strategy of human resources management;
- Set the mode of blocking (sanctioning) the fulfillment of a certain set of functional duties of personnel from job descriptions;
- Estimation the integral state of the parameters characterizing the basic strategy of human resources management (hereinafter referred to as HR strategy) and HR with the influence of sanctions.

The following terms and abbreviations are used in the work: strategy of human resources management - a set of functional duties of the enterprise personnel, fixed in the job description;

To analyze the state of an economic object (enterprise) by the method of integrated indicators, we represent the economic object as a system $S$ and give its description.

## 2. Method

It is enough to imagine the economic activity of the enterprise as $S=\{T,X\}$, where $T=\{t:t=1,...,T_{max}\}$ - a lot of time points with a selected interval for analysis; $X$ - space of system parameters; $x(t)=[x^1(t),x^2(t),...,x^n(t)]^T \in X$ – $n$ – vector of indicators corresponding to the state of the system. Indicators of the vector $x^i(t)$ - the value of financial expenses and income of the enterprise. The dimension of the system $n$ is 1.2 million parameters. Based on the parameters $X$ and $T$, we consider our economic object a multidimensional dynamic system (hereinafter referred to as the system).

The analysis of the system at the moment $t$ is performed $x(t)$ on the basis $k$ of previous measures. The parameter $k$ is the length of the time series segment (accepted $k=6$ months in the work). Then we have a matrix

$$X_k(t) = \begin{bmatrix} x^T(t-1) \\ x^T(t-2) \\ \dots \\ x^T(t-k) \end{bmatrix} = \begin{bmatrix} x^1(t-k) & x^2(t-k) & \cdots & x^n(t-k) \\ x^1(t-k) & x^2(t-k) & \cdots & x^n(t-k) \\ \dots & \dots & \dots & \dots \\ x^1(t-k) & x^2(t-k) & \cdots & x^n(t-k) \end{bmatrix} \quad (1)$$

$$R_k(t) = \frac{1}{k-1} \overset{o}{X}^T_k(t) \overset{o}{X}_k(t) = \left\| r_{ij}(t) \right\|, \quad (2)$$

$$r_{ij}(t) = \frac{1}{k-1} \sum_{l=1}^{k} \overset{o}{x}^i(t-l) \overset{o}{x}^j(t-l), \quad i,j=1,...,n, \quad (3)$$

where $t$ are the time instants, $r_{ij}(t)$ are the correlation coefficients of the variables $x^i(t)$ и $x^j(t)$ at the time instant $t$.

Next we form one of the four integral indicators – the sum of the absolute indicators of the correlation coefficients. It is indicator for express estimation of the correlation of system parameters $G_i(t)$:

$$R_i(t) = G_i(t) = \sum_{j=1}^{n} |r_{ij}(t)|. \tag{4}$$

The state of the entire system is calculated as:

$$G = \sum_{t=1}^{T=\max} \sum_{i=1}^{n} G_i(t). \tag{5}$$

The strategy of human resources management $V$ is a set of functional instructions for all personnel, which can be represented $V_i^k$ as a set of personnel management strategies:

$$V = \sum_{t=1}^{T=\max} \sum_{i=1}^{n} V_i^k(t). \tag{6}$$

We will carry out the identification of the performed functions of the system with the fulfilled functional duties of the personnel, then each strategy of human resources management $V_i^k$ is characterized by the functional responsibilities of the personnel:

$$V_i^k = \sum_{i=1}^{n} v_i^j(x_j^i) \to \min, \tag{7}$$

where $v_i^j$ is the functional duty of the employee prescribed in the job description (compliance $x_j^i$ is $v_i^j$ set as 1-yes, 0-no); $x_j^i$ - the costs of the $i$ - functional obligation at the enterprise for the employee $j$. The level of automation of functional responsibilities is taken into account as the work performed as part of the job description of the operator of the automated process.

Then, sanctions are a set of blocked functional responsibilities in a vector of values $v(t) = \left[v^1(t), v^2(t), ..., v^n(t)\right]^T \in V$ - $n$ - a dimension characterizing the state of the HR.

Payment of the functional duties of employees of the economic system is limited by resources $C$, then $C(X) \leq C$. This restriction applies to all subsystems of the researched system.

The implementation of the method is performed in the author's complex of programs.

## 3. Characteristics of the research objects

In the first example, a construction company established in 2002 is used as an economic object for research. The company carried out the construction of 4 residential buildings. It was monitored from December 2003 until mid-2008 (52 periods). The need to use this enterprise, as an example, is due to the fact that it is necessary to record changes in the state of the system from changes in the parameters characterizing the result of applying control actions to the functional responsibilities of personnel (strategy of human resources management) in a well-researched, real-life economic facility. Simply put, already accepted control actions are not reviewed so as not to mix the system response to them with the studied control actions (functional responsibilities of the staff). At this facility, the following have already been determined: the boundaries of the system space, the matrices that determine the structure of the data space and the space of control actions, the possibility of using an integral indicator to assess the state of an economic object and its quality assessment of control actions. The space of the system is characterized by the following parameters: description of the project and owners, fixed assets, product characteristics, market analysis, organizational plan, marketing data, production structure (warehouses, logistics), analysis of resource consumption rates, environmental impact, analysis of project risks, financial income model and expenses. Important events in the system in figure 1: approval of working documentation in the 22nd period, obtaining a building license in the 32nd period, the implementation of a quality management system in the 34th period, the opening of a credit line in the 40th period, the start of construction in the 42nd period, the opening of a subsidiary in 44 period. Since the 34th period the company has been affected by the global financial crisis. In 2006, the author personally implemented an automated data control system at the enterprise to analyze its activities. The above parameters are entered into the software package [20, 21], to simulate the state of the system at different points in time with different control actions.

In the second example, an investment project of the Russian Federation on deep wood processing is considered. The research enterprise is included in the list of priority projects of Russia. Since 2007 priority projects have been formed. It includes the most technologically advanced enterprises and projects in the forest industry. Such projects undergo a thorough expertise of the ministries in the regional and federal levels. A status of priority project of the Russian Federation gives the enterprise tax benefits. This procedure is governed by the Resolution of the Russian Federation of June 30, 2007 No. 419 "On priority investment projects in the field of forest development"

The investigated enterprise prepares 800 thousand cubic meters of industrial roundwood in the North of the Krasnoyarsk Territory (upper warehouse). Further, along the Yenisei River, industrial roundwood is delivered for deep processing and production of wood products. This type of deep wood processing is waste-free.

Industrial roundwood stocks are formed in such a way as to ensure uninterrupted operation of the enterprise for 1.5 years.

## 4. Experiment result

Example 1. Without a sanctions regime, the system has an initial set of fulfilled functional duties of personnel characterizing the HR ($V_1^6$ - the basic strategy) and the value of the indicator is 153,080 in figure 1. A restrictive regime (sanctions) has been introduced through the exclusion of engineering personnel from activities important for construction. The list of functional duties of personnel blocked by sanctions: engineer-designer in the period 1-19 - providing technical documentation for the construction site and in periods 38-42, 50-54, 62-66 - development of technical specifications; engineer of the technical department in the period 1-20 - development of technical solutions, in the period 28-36 - development of projects for the production of works; design engineer in the period 1-27 - participation in the installation of structures, tests and commissioning, in the period 28-70, in the period 71-73 - participation in the installation of structures. In table 1 the entire list of functional responsibilities blocked by sanctions is presented

**Table 1.** Sanctioned functional responsibilities strategy ($V_1^6$).

| Period (t) | Employee's position | Functional (staff responsibilities) $v_i^j$ |
|---|---|---|
| 1 - 19 | Concept engineer | Providing technical documentation for the construction site |
| 38 - 42 | Concept engineer | Development of technical specifications |
| 50 - 54 | Concept engineer | Development of technical specifications |
| 62 - 66 | Concept engineer | Development of technical specifications |
| 1 - 20 | Proof engineer | Creation of technical solutions |
| 28 - 36 | Proof engineer | Creation of work production projects |
| 1 - 27 | Facility design engineer | Control of installation of structures, testing and commissioning of facilities |
| 28 - 70 | Facility design engineer | Designing the main sections of the project |
| 71 - 73 | Facility design engineer | Control and participation in the installation of structures |
| 1 - 32 | Comprehensive services | Implementation of construction supervision, participation in seminars and conferences |
| 50 – 51 | Comprehensive services | Author's supervision |
| 62 - 63 | Comprehensive | Author's supervision |
| | | services |
| 1 – 18 | Project chief engineer | Project Coordination |
| 66 - 73 | Project chief engineer | Project Coordination |
| 38 - 39 | Project chief engineer | Author's supervision |
| 1 - 20 | Section foremaster | Implementation of technical control over the implementation of construction and installation work; acceptance of completed volumes |
| 25 – 26 | Section foremaster | Material Accounting |
| 38 – 39 | Section foremaster | Material Accounting |
| 50 – 51 | Section foremaster | Material Accounting |
| 62 - 63 | Section foremaster | Material Accounting |
| 1 - 20, 25, 38, 50, 62 | Site supervisor | Performing construction and installation work |
| 1 - 20 | Production and Technical Department Engineer | Determination of the volume of work performed |
| 1 - 20 | Proof engineer | Keeping records of completed construction and assembly works and reporting on the implementation of plans |
| 1 - 20 | Occupational safety engineer | Monitoring the implementation of labor protection measures |

A figure 1 shows blocked functional responsibilities of engineering personnel $v_i^j$ affect the performance $x_j^i$ of system functions (7). The value of the HR in the sanctions regime ($V_2^6$ - Strategy 2) is 155,896.

Functions blocked by sanctions can be restored within 2 months by attracting the services of third-party organizations in the amount of 148.5 million rubles, which is equivalent to a 10% rise in the financial cost of the entire project for 4 years.

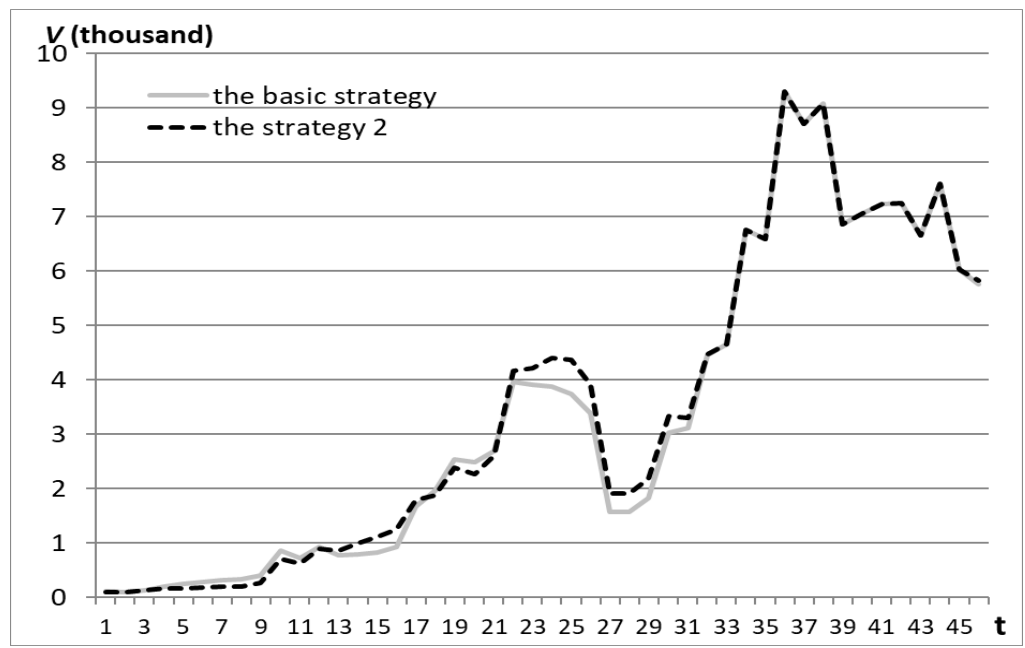

**Figure 1.** Parameter dynamics of $V_i(t)$.

Then the influence of sanctions on the system is measured as the difference between the state of the control system parameters and is equal to minus 2,816. The value is due to the higher multi-factorial response of the system to more difficult conditions for it. The multifactorial effect in systems is similarly disclosed in a separate paper [18]. The loss of the enterprise from the imposed sanctions over 4 years will amount to 61.63 million rubles.

Example 2. In figure 2 the normal mode of implementation of the investment project is characterized by strategy 1 under number 1. A restrictive regime is introduced in the form of sanctions on the provision of personnel to complete the installation of deep wood processing equipment. To start the investment project, the Strategy 2 is applied, which involves the search and involvement of other foreign employees for the installation of equipment. The table 2 shows different between of the search foreign employees for the installation of equipment from basic strategy.

**Table 2.** The impact of sanctions $\Delta V$ on the economic parameters of the investment project.

| № | Parameters | $\Delta V$ | Strategy | |
|---|---|---|---|---|
| | | | $V_1^6$ | $V_2^6$ |
| 1 | Line of credit | +26% | 100% | 126% |
| 2 | Loan rate | +3% | 10% | 13% |
| 3 | Owner's investments | - | - | - |
| 4 | State benefits and subsidies | +67% | 33% | 100% |
| 5 | Sale of industrial roundwood (months) | +5 | From 5 | From 10 |
| 6 | Product sales (months) | +6 | From 21 | From 27 |
| 7 | Fixed asset delivery offset (months) | +5 | - | From 5 |

## 5. The discussion of the results

The result of the experiment is one of blocking the performance of staff functional duties is conditionally identical to the sanctions imposed by Austria against one of Gazprom's subsidiaries. According to public sources from 10.16.2019, it is known that the Austrian company LMF forcibly switched off the compressor via satellite, thereby blocking the functional responsibilities of the engineering personnel for pumping gas. Gazprom said it has the ability to replace equipment and will solve this problem no later than 5 months. The result of experiment two is conditionally identical to the implementation of the Nord Stream 2 project under the conditions of sanctions.

Another real process, which is conditionally identical to the sanctions introduced in the experiment, is the emigration of engineering and technical personnel together with the family from the country to another job, which hinders the development of the economic object for 2-4 months. The damage to economic activity is caused not only to the enterprise where there is a halt to the fulfillment of functional duties, but also to a larger economic system: region, territory, country. This synergistic damage has been little researched, which may serve as one of the incentives for further research in this direction. This work was completed as part of a series of articles on the control of economic objects [22-27] and the sanctions regime for them [28, 29].

## 4. Conclusion

The tasks set at the beginning of the work are completed:
- An economic object (enterprise) was set as a multidimensional dynamic system with respect to the parameters of set *X* and *T*.

- Assigned to the performed functions of the system functional $x_j^i$ duties from job descriptions $v_i^j$, thereby forming a set of job descriptions characterizing the basic strategy of human resources management $V_i^k$. The HR was defined as $V_i^k = \sum_{i=1}^{n} v_i^j(x_j^i)$.
- The mode of blocking (sanctions) of the performance of functional duties $v_i^j$ from the job descriptions of the personnel was set.
- An integral indicator was evaluated for two states, the basic HR $V_1^6$ - 153,080 and the HR with sanctions $V_2^6$ - 155,896 described by the vector $v(t) = \left[ v^1(t), v^2(t), ..., v^n(t) \right]^T \in V$.

The experimental data are compared with real situations.

The goal set at the beginning of the work, to assess the state of an economic object by an integral indicator as a multidimensional dynamic system, was achieved when the staff performed their functional duties in the base mode and sanctions mode with unknown environmental parameters.